\documentclass[conference, 9pt]{IEEEtran}
\IEEEoverridecommandlockouts
\usepackage{graphicx}
\usepackage{color}
\usepackage{placeins}
\usepackage{float}
\usepackage{tabularx,colortbl}
\usepackage{amsmath, amssymb}
\usepackage{epsfig}

\usepackage{amsthm}

\begin{document}
\title{Conservative Signal Processing Architectures\\For Asynchronous, Distributed Optimization\\Part I: General Framework}

\author{\IEEEauthorblockN{Thomas A. Baran and Tarek A. Lahlou\thanks{The authors wish to thank Analog Devices, Bose Corporation, and Texas Instruments for their support of innovative research at MIT and within the Digital Signal Processing Group.}}
\IEEEauthorblockA{Digital Signal Processing Group\\
Massachusetts Institute of Technology }}

\maketitle

\begin{abstract}
This paper presents a framework for designing a class of distributed, asynchronous optimization algorithms, realized as signal processing architectures utilizing various conservation principles.  The architectures are specifically based on stationarity conditions pertaining to primal and dual variables in a class of generally nonconvex optimization problems.  The stationarity conditions, which are closely related to the principles of stationary content and co-content that can be derived using Tellegen's theorem in electrical networks, are in particular transformed via a linear change of coordinates to obtain a set of linear and nonlinear maps that form the basis for implementation.  The resulting algorithms specifically operate by processing a linear superposition of primal and dual decision variables using the associated maps, coupled using synchronous or asynchronous delay elements to form a distributed system.  A table is provided containing specific example elements that can be assembled to form various optimization algorithms directly from the corresponding problem statements.

\end{abstract}
\begin{IEEEkeywords} Asynchronous optimization, distributed optimization, conservation \end{IEEEkeywords}

\IEEEpeerreviewmaketitle

\vspace{-1mm}
\section{Introduction}

In designing distributed, asynchronous algorithms for optimization, a common approach is to begin with a non-distributed iteration or with a distributed, synchronous implementation and attempt to organize variables so that the algorithm distributes across multiple unsynchronized processing nodes.\cite{NealBoyd}\cite{WeiOzdaglar1}\cite{ForeroCanoGiannakis}  An important limitation of this research strategy is that it does not generally involve any particular systematic approach for performing such an organization.  The presented framework addresses this by introducing techniques for directly designing a variety of algorithm architectures for convex and nonconvex optimization that naturally distribute across multiple processing elements utilizing synchronous or asynchronous updates.

This paper is one of two parts.  In particular this paper establishes the general framework and provides a straightforward strategy for designing distributed, asynchronous optimization algorithms directly from associated problem statements.  Part II \cite{BaranLahlouPartII} provides examples of this strategy, a discussion of convergence, as well as simulations of various resulting algorithms.

\subsection{Classes of maps}

Following the convention suggested in \cite{Willems1}, we make use of several specific terms in describing linear and nonlinear maps.  The term ``neutral'' will refer to any map $m(\cdot)$ for which
\begin{equation}
\label{eq:neturalDef}
||m(\underbar{x})|| = ||\underbar{x}||,\ \ \forall \underbar{x},
\end{equation}
with $||\cdot ||$ being used here and throughout this paper to denote the $2$-norm. The expression ``$\forall \underbar{x}$'' in Eq.~\ref{eq:neturalDef} is used to indicate all vectors $\underbar{x}$ in the domain over which $m(\cdot)$ is defined.

We will denote as ``passive about $\underbar{x}'$'' any map $m(\cdot)$ for which
\begin{equation}
\label{eq:passiveDef}
\sup_{\underbar{x}\neq 0} \frac{||m(\underbar{x}+\underbar{x}') - m(\underbar{x}')||}{||\underbar{x}||} \leq 1.
\end{equation}
As a subset of passive maps, we will denote as ``dissipative about $\underbar{x}'$'' any map $m(\cdot)$ for which
\begin{equation}
\label{eq:dissipativeDef}
\sup_{\underbar{x}\neq 0} \frac{||m(\underbar{x}+\underbar{x}') - m(\underbar{x}')||}{||\underbar{x}||} < 1.
\end{equation}
A map that is ``passive everywhere'' or ``dissipative everywhere'' is a map that is passive, or respectively dissipative, about all points $\underbar{x}'$.

The term ``source'' will be used to refer to a map that is written as
\begin{equation}
\label{eq:sourceEq}
m(\underbar{d}) = S \underbar{d} + \underbar{e},
\end{equation}
where $\underbar{e}$ is a constant vector and where the map that is associated with the matrix $S$ is passive.

\subsection{Notation for partitioning vectors}
\label{ssec:notation}

We will commonly refer to various partitionings of column vectors, each containing a total of $N$ real scalars, in the development and analysis of the presented class of architectures.  To facilitate the indexing associated with this, we establish an associated notational convention.  Specifically we will refer to two key partitionings of a length-$N$ column vector $\underbar{z}$, indicated using superscripts whose meanings will be discussed in Section \ref{sec:classArch}. In one such partitioning the elements are arranged into a total of $K$ column vectors denoted $\underbar{z}_k^{(CR)}$, and in the other the elements are partitioned into a total of $L$ column vectors denoted $\underbar{z}_{\ell}^{(LI)}$.  Each vector $\underbar{z}_{\ell}^{(LI)}$ will also be partitioned into subvectors denoted $\underbar{z}^{(i)}_{\ell}$ and $\underbar{z}^{(o)}_{\ell}$.  We write all of this formally as
\begin{align}
[z_1, \dots, z_N ]^T  &= [\underbar{z}^{(CR)^T}_1, \dots, \underbar{z}^{(CR)^T}_K ]^T \label{eq:partSchemeFirst}\\
  &= [\underbar{z}^{(LI)^T}_1, \dots, \underbar{z}^{(LI)^T}_L ]^T \\
  &= \underbar{z} \in \mathbb{R}^N .
\end{align}
\begin{equation}
z_{\ell}^{(LI)} = [\underbar{z}^{(i)^T}_{\ell}, \underbar{z}^{(o)^T}_{\ell}]^T,\ \ \ell = 1,\dots L.
\end{equation}
The length of a particular subvector $\underbar{z}_k^{(CR)}$, $\underbar{z}_{\ell}^{(LI)}$, $\underbar{z}^{(i)}_{\ell}$, or $\underbar{z}^{(o)}_{\ell}$ will respectively be denoted $N_k^{(CR)}$, $N_{\ell}^{(LI)}$, $N_{\ell}^{(i)}$, $N_{\ell}^{(o)}$, with
\begin{align}
N &= N_1^{(CR)} + \dots + N_K^{(CR)}\\
 &= N_1^{(LI)} + \dots + N_L^{(LI)}
\end{align}
\begin{equation}
N_{\ell}^{(LI)} = N_{\ell}^{(i)} + N_{\ell}^{(o)},\ \ \ell = 1,\dots L. \label{eq:partSchemeLast}
\end{equation}

\section{Class of optimization problems}

The class of optimization problems addressed within the presented framework is similar in form to those problems described by the well-known principles of stationary content and co-content in electrical networks,\cite{Millar}\cite{PenfieldSpenceDuinker} which have been used in constructing circuits for performing convex and nonconvex optimization.\cite{ChuaLin}\cite{Dennis}\cite{KennedyChua}\cite{Wyatt1995}  These principles and implementations implicitly or explicitly utilize a nonconvex duality theory where physical conjugate variables, e.g.~voltage and current, are identified as primal and dual decision variables within the associated network.  In this paper we will specifically utilize the multidimensional, parametric generalization of the principles of stationary content and co-content that was developed in \cite{tbaran-phd}. 

We define a dual pair of problems within the presented class first in a form that will be used for analysis from a variational perspective, which we will refer to as ``canonical form''.  We will also utilize an alternative form obtained by performing algebraic manipulations on problems in canonical form, referred to as ``reduced form''.  Optimization problems will typically be written in reduced form for the purpose of relating their formulations to those of generally well-known classes of convex and nonconvex problems.

\subsection{Canonical-form representation}

Making use of the partitioning convention established in Eqns.~\ref{eq:partSchemeFirst}-\ref{eq:partSchemeLast}, we write a specific primal problem in canonical form as
\begin{eqnarray}
\label{eq:optim-primal}
 &\displaystyle \min_{ \{ y_1,\dots, y_N \} \atop \{ a_1,\dots, a_N \} } & \sum_{k=1}^K  Q_k(\underbar{y}^{(CR)}_k) \label{eq:optim-primal1}\\
 &\mbox{s.t.} & \underbar{a}^{(CR)}_k = f_k(\underbar{y}^{(CR)}_k),\ \ k=1,\dots, K \label{eq:optim-primal2} \\
 &\ & A_{\ell} \underbar{a}^{(i)}_{\ell} = \underbar{a}^{(o)}_{\ell},\ \ \ell=1,\dots, L \label{eq:optim-primal3}.
\end{eqnarray}
The functionals $Q_k(\cdot): \mathbb{R}^{N^{(CR)}_k} \to \mathbb{R}$ composing the summation in (\ref{eq:optim-primal1}) are in particular related to the functions $f_k(\cdot): \mathbb{R}^{N^{(CR)}_k} \to \mathbb{R}^{N^{(CR)}_k}$ in (\ref{eq:optim-primal2}) according to the following:
\begin{equation}
\label{eq:contentQual}
\nabla Q_k(\underbar{y}^{(CR)}_k) = J_{f_k}^T(\underbar{y}^{(CR)}_k)g_k(\underbar{y}^{(CR)}_k),
\end{equation}
where $f_k(\cdot)$ and $g_k(\cdot): \mathbb{R}^{N^{(CR)}_k} \to \mathbb{R}^{N^{(CR)}_k}$ are generally nonlinear maps whose respective Jacobian matrices $J_{f_k}(\underbar{y}^{(CR)}_k)$ and $J_{g_k}(\underbar{y}^{(CR)}_k)$ are assumed to exist.\footnote{We use the convention that the entry in row $i$ and column $j$ of $J_{f_k}(\underbar{y}^{(CR)}_k)$ is the partial derivative of output element $i$ of $f_k(\underbar{y}^{(CR)}_k) $, with respect to element $j$ of the input vector $ \underbar{y}^{(CR)}_k $, evaluated at $\underbar{y}^{(CR)}_k$.} Each of $A_{\ell}: \mathbb{R}^{N^{(i)}_{\ell}} \to \mathbb{R}^{N^{(o)}_{\ell}} $, $\ell = 1, \dots, L$, is a linear map.

Given a primal problem written in canonical form as (\ref{eq:optim-primal1})-(\ref{eq:optim-primal3}), we write the associated dual problem in canonical form as
\begin{eqnarray}
\label{eq:optim-dual}
 &\displaystyle \max_{\{ y_1,\dots, y_N \} \atop \{ b_1,\dots, b_N \} } & -\sum_{k=1}^K R_k(\underbar{y}^{(CR)}_k) \label{eq:optim-dual1}\\
 &\mbox{s.t.} & \underbar{b}_k = g_k(\underbar{y}^{(CR)}_k),\ \ k=1,\dots, K \label{eq:optim-dual2} \\
 &\ & \underbar{b}^{(i)}_{\ell} = -A^T_{\ell}  \underbar{b}^{(o)}_{\ell},\ \ \ell=1,\dots, L, \label{eq:optim-dual3}
\end{eqnarray}
where
\begin{equation}
\label{eq:contentCoContentQual}R_k(\underbar{y}^{(CR)}_k) = \left\langle f_k(\underbar{y}^{(CR)}_k),g_k(\underbar{y}^{(CR)}_k)\right\rangle -Q_k(\underbar{y}^{(CR)}_k),\ \ k=1,\dots, K,
\end{equation}
and with $\langle \cdot , \cdot \rangle$ denoting the standard inner product.  As is suggested by the notation established in Subsection \ref{ssec:notation}, the primal and dual costs and constraints in (\ref{eq:optim-primal1}), (\ref{eq:optim-primal2}), (\ref{eq:optim-dual1}), and (\ref{eq:optim-dual2}) will be specified using a total of $K$ constitutive relations within the presented class of architectures.  Likewise the primal and dual linear constraints in (\ref{eq:optim-primal3}) and (\ref{eq:optim-dual3}) will be specified in the presented class of architectures using a total of $L$ linear interconnection elements.

\subsection{Reduced-form representation}

For various choices of $Q_k(\cdot)$ and $f_k(\cdot)$, it is generally possible that the set of points traced out in $\underbar{a}^{(CR)}_k$-$Q_k$, generated by sweeping $\underbar{y}^{(CR)}_k$, is one that could equivalently have been generated using a functional relationship mapping from $\underbar{a}^{(CR)}_k\in \mathbb{R}^{N^{(CR)}_k}$ to $Q_k\in \mathbb{R}$, possibly with $\underbar{a}^{(CR)}_k$ being restricted to an interval or set.  In cases where this is possible for all $f_k$-$Q_k$ pairs forming (\ref{eq:optim-primal1})-(\ref{eq:optim-primal3}), we will formulate the problem in terms of functionals $\widehat{Q}_k(\cdot): \mathbb{R}^{N^{(CR)}_k} \to \mathbb{R}$ and sets $\mathcal{A}_k\subseteq \mathbb{R}^{N^{(CR)}_k}$ in what we refer to as ``reduced form'':
\begin{eqnarray}
\label{eq:optim-primalRed1}
 &\displaystyle \min_{\{a_1,\dots, a_N\} } & \sum_{k=1}^K \widehat{Q}_k(\underbar{a}^{(CR)}_k)\\
 &\mbox{s.t.} & \underbar{a}^{(CR)}_k \in \mathcal{A}_k,\ \ k=1,\dots, K \\
 &\ & A_{\ell} \underbar{a}^{(i)}_{\ell} = \underbar{a}^{(o)}_{\ell},\ \ \ell=1,\dots, L \label{eq:optim-primalRed3}.
\end{eqnarray}
A reduced-form representation may specifically be used when $Q_k(\cdot)$, $f_k(\cdot)$, $\widehat{Q}_k(\cdot)$, and $\mathcal{A}_k$ satisfy the following relationship:
\begin{equation}
\left\{ \left[ \begin{array}{c} f_k(\underbar{y}^{(CR)}_k) \\ Q_k(\underbar{y}^{(CR)}_k) \end{array} \right] : \underbar{y}^{(CR)}_k\in \mathbb{R}^{N^{(CR)}_k} \right\} = \left\{ \left[ \begin{array}{c} \underbar{a}^{(CR)}_k \\ \widehat{Q}_k(\underbar{a}^{(CR)}_k) \end{array} \right] : \underbar{a}^{(CR)}_k\in \mathcal{A}_k \right\}. \label{eq:reducedFormPrimalEquivalence}
\end{equation}

The key idea in writing a problem in reduced form, i.e.~(\ref{eq:optim-primalRed1})-(\ref{eq:optim-primalRed3}), is to provide a formulation that allows for set-based constraints on decision variables, in addition to allowing for cost functions that need not be differentiable everywhere.  It is, for example, generally possible to define functions $f_k(\cdot)$ and $g_k(\cdot)$ that are differentiable everywhere, resulting in a canonical-form cost term $Q_k(\cdot)$ that is differentiable everywhere, and for an associated reduced-form cost term $\widehat{Q}_k(\cdot)$ satisfying Eq.~\ref{eq:reducedFormPrimalEquivalence} to have knee points where its derivative is not well-defined.  This issue is discussed in greater detail in \cite{tbaran-phd}.

A dual canonical-form representation (\ref{eq:optim-dual1})-(\ref{eq:optim-dual3}) may similarly be written in reduced form:
\begin{eqnarray}
 &\displaystyle \max_{\{ b_1,\dots, b_N \} } & -\sum_{k=1}^K \widehat{R}_k(\underbar{b}_k) \label{eq:optim-dualRed1}\\
 &\mbox{s.t.} & \underbar{b}_k \in \mathcal{B}_k,\ \ k=1,\dots, K \label{eq:optim-dualRed2} \\
 &\ & \underbar{b}^{(i)}_{\ell} = -A^T_{\ell}  \underbar{b}^{(o)}_{\ell},\ \ \ell=1,\dots, L, \label{eq:optim-dualRed3}
\end{eqnarray}
where $\widehat{R}_k(\cdot): \mathbb{R}^{N^{(CR)}_k} \to \mathbb{R}$ and $\mathcal{B}_k\subseteq \mathbb{R}^{N^{(CR)}_k}$ for which
\begin{equation}
\left\{ \left[ \begin{array}{c} g_k(\underbar{y}^{(CR)}_k) \\ R_k(\underbar{y}^{(CR)}_k) \end{array} \right] : \underbar{y}^{(CR)}_k\in \mathbb{R}^{N^{(CR)}_k} \right\} = \left\{ \left[ \begin{array}{c} \underbar{b}_k \\ \widehat{R}_k(\underbar{a}^{(CR)}_k) \end{array} \right] : \underbar{b}_k\in \mathcal{B}_k \right\}. \label{eq:reducedFormDualEquivalence}
\end{equation}
We note that if a primal problem is representable in reduced form, the dual problem may or may not have an associated reduced-form representation, or vice-versa.  The last row of the table in Fig.~\ref{fig:bigTable} provides an example of this.

\subsection{Stationarity conditions}

As a consequence of the formulation of the primal and dual problems in canonical form, respectively (\ref{eq:optim-primal1})-(\ref{eq:optim-primal3}) with (\ref{eq:contentQual}), and (\ref{eq:optim-dual1})-(\ref{eq:optim-dual3}) with (\ref{eq:contentCoContentQual}), the dual pair of feasibility conditions serve as stationarity conditions for the dual pair of costs.  Specifically, any point described by the set of vectors ${\underbar{y}^{\star}_k}^{(CR)}$ that satisfies Eqns.~\ref{eq:optim-primal2}-\ref{eq:optim-primal3} and \ref{eq:optim-dual2}-\ref{eq:optim-dual3}, is a point about which both the primal cost (\ref{eq:optim-primal1}) and dual cost (\ref{eq:optim-dual1}) are constant to first order, given any small change in ${\underbar{y}^{\star}_k}^{(CR)}$ for which the primal constraints (\ref{eq:optim-primal3}) and dual constraints (\ref{eq:optim-dual3}) remain satisfied.  A proof of essentially this statement, which is a multidimensional generalization of the well-known principles of stationary content and co-content in electrical networks \cite{Millar,PenfieldSpenceDuinker}, can be found in \cite{tbaran-phd}.

\vspace{-1mm}
\section{Class of architectures}
\label{sec:classArch}

The key idea behind the presented class of architectures is to determine a solution to the stationarity conditions composed of Eqns.~\ref{eq:optim-primal2}-\ref{eq:optim-primal3} and \ref{eq:optim-dual2}-\ref{eq:optim-dual3}, in particular by interconnecting various signal-flow elements and running the interconnected system until it nears a fixed point.  The elements in the architecture are specifically memoryless, generally nonlinear maps that are coupled via synchronous or asynchronous delays, which we will model as discrete-time,  sample-and-hold elements triggered in the asynchronous case by independent discrete-time Bernoulli processes.

\vspace{-1mm}
\begin{figure}[h!]
\begin{center}
\epsfig{file=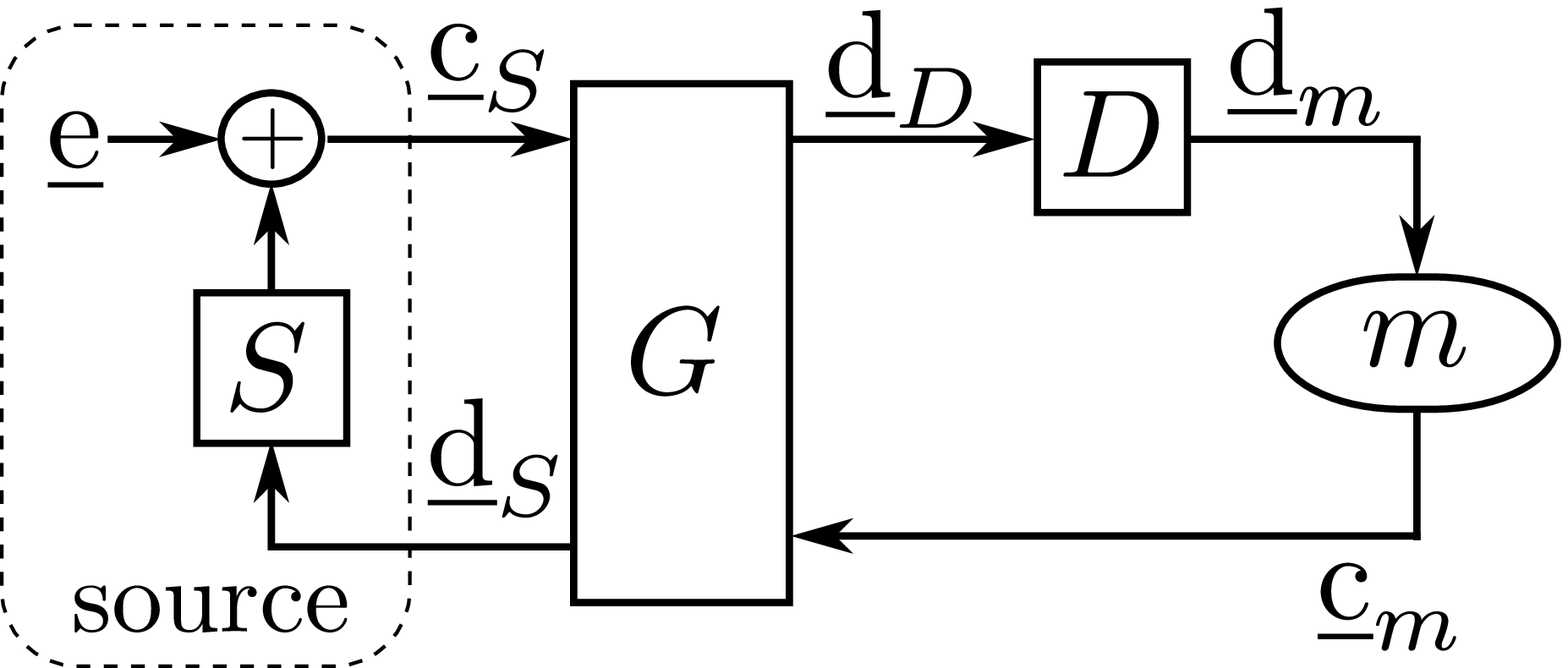,width=1.5in}
\vspace{-1mm}
\caption{General interconnection of elements in the presented architectures.\label{fig:convergence}}
\end{center}
\end{figure}

\vspace{-1mm}
The approach for interconnecting the various system elements is depicted in Fig.~\ref{fig:convergence}.  Referring to this figure, systems in the presented class of architectures will be composed of a set of $L$ memoryless, neutral, linear interconnections (LI) denoted $G_{\ell}$ and in the aggregate denoted $G$, coupled directly to a set of $K$ maps $m_k(\cdot)$, referred to as constitutive relations (CRs).  A subset of the maps $m_k(\cdot)$ that have the property of being source elements are specifically connected directly to $G$, and the remaining maps $m_k(\cdot)$, denoted on the whole as $m(\cdot)$, are coupled to the interconnection via delay elements.  Algebraic loops will generally exist between the remaining source elements and the interconnection, and as these are linear may be eliminated by performing appropriate algebraic reduction.

Given a particular system within the presented class, we have two key requirements of the system:
\begin{itemize}
\item[(R1)] The system converges to a fixed point, and
\item[(R2)] Any fixed point of the system corresponds to a solution of the stationarity conditions in Eqns.~\ref{eq:optim-primal2}-\ref{eq:optim-primal3} and \ref{eq:optim-dual2}-\ref{eq:optim-dual3}.
\end{itemize}
The issue of convergence in (R1) relates to the dynamics of the interconnected elements, and (R2) relates to the behavior\footnote{Consistent with the convention in \cite{Willems2}, we refer to the ``behavior'' of a system of maps as the set of all input and output signal values consistent with the constraints imposed by the system.  The term ``graph form'' has also been used to denote a similar concept.\cite{NealBoyd}} of the interconnection of the various memoryless maps composing the system, with the delay elements being replaced by direct sharing of variables.

\begin{figure*}[t!]
\begin{center}
\vspace{-4mm}
\epsfig{file=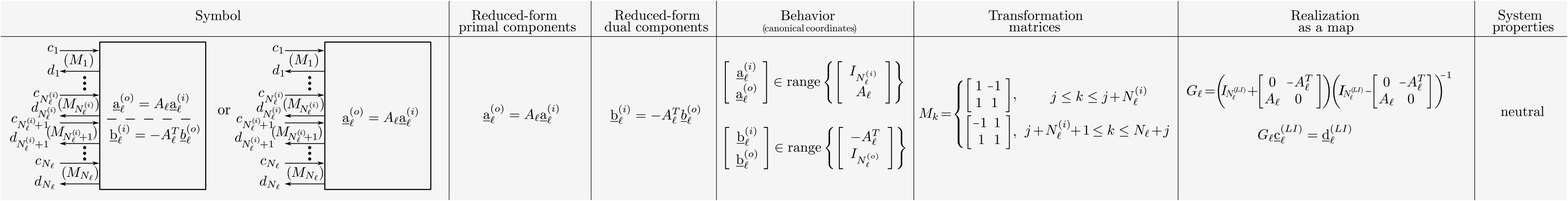,width=7in}
\vspace{-0.1in}
\caption{Example LI elements, graphically denoted using rectangles, satisfying Eq.~\ref{eq:Atransformation}.  The maps in column 6 are used in implementation. \label{fig:smallTable}}
\vspace{0.13in}
\epsfig{file=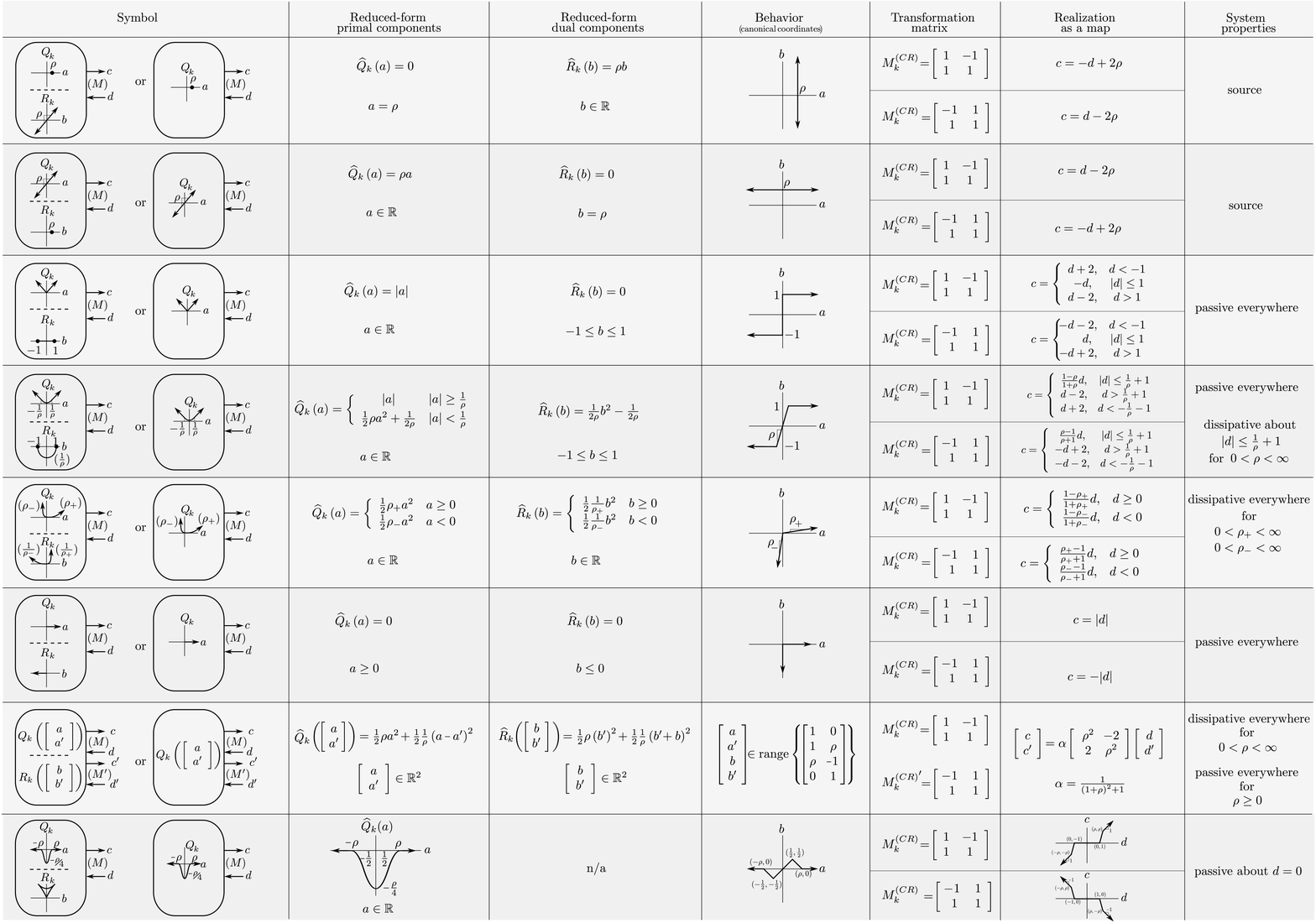,width=7in}
\vspace{-0.1in}
\caption{Example CR elements, graphically denoted using rounded rectangles, satisfying Eq.~\ref{eq:fgTransformation}.  The maps in column 6 are used in implementation.\label{fig:bigTable}}
\end{center}
\vspace{-6mm}
\end{figure*}

\vspace{-1mm}
\subsection{Coordinate transformations}

In satisfying (R1) and (R2), the general strategy is to perform a linear, invertible coordinate transformation of the primal and dual decision variables $\underbar{a}$ and $\underbar{b}$, and to use the transformed stationarity conditions, obtained by transforming Eqns.~\ref{eq:optim-primal2}-\ref{eq:optim-primal3} and \ref{eq:optim-dual2}-\ref{eq:optim-dual3}, to form the basis for the synchronous or asynchronous system summarized in Fig.~\ref{fig:convergence}.  The linear stationarity conditions in Eqns.~\ref{eq:optim-primal3} and \ref{eq:optim-dual3} will in particular be used in defining the linear interconnections $G_k$, and the generally nonlinear stationarity conditions in Eqns.~\ref{eq:optim-primal2} and \ref{eq:optim-dual2} will be used in defining the constitutive relations $m_k(\cdot)$. 

We specifically utilize coordinate transformations consisting of a pairwise superposition of the primal and dual decision variables $a_i$ and $b_i$, resulting in transformed variables denoted $c_i$ and $d_i$.  The associated change of coordinates is written formally in terms of a total of $N$, $2\times 2$ matrices $M_i$ as
\begin{equation}
\left[ \begin{array}{c} c_i \\ d_i \end{array} \right] = M_i \left[ \begin{array}{c} a_i \\ b_i \end{array} \right],\ \ i = 1,\dots, N. \label{eq:linTransformation}
\end{equation}
Viewing the transformed variables $c_i$ and $d_i$ as entries of column vectors written $\underbar{c}$ and $\underbar{d}$, we will make use of the partitioning scheme described in Eqns.~\ref{eq:partSchemeFirst}-\ref{eq:partSchemeLast}.  Linear maps denoted $M_k^{(CR)}$ and 
$M_{\ell}^{(LI)}$ will likewise be used to represent the relationship described in Eq.~\ref{eq:linTransformation} in a way that is consistent with the various associated partitionings:
\begin{align}
\left[ \begin{array}{c} \underbar{c}^{(CR)}_k \\ \underbar{d}^{(CR)}_k \end{array} \right] &= M^{(CR)}_k \left[ \begin{array}{c} \underbar{a}^{(CR)}_k \\ \underbar{b}^{(CR)}_k \end{array} \right],\ \ k = 1,\dots, K\\
\left[ \begin{array}{c} \underbar{c}^{(LI)}_{\ell} \\ \underbar{d}^{(LI)}_{\ell} \end{array} \right] &= M^{(LI)}_{\ell} \left[ \begin{array}{c} \underbar{a}^{(LI)}_{\ell} \\ \underbar{b}^{(LI)}_{\ell} \end{array} \right],\ \ \ell = 1,\dots, L.
\end{align}

Referring to Fig.~\ref{fig:convergence}, we will use the variables $c_i$ and $d_i$ to respectively denote the associated linear interconnection inputs and outputs, and we will denote the constitutive relation inputs using $d_i$ and the associated outputs using $c_i$.  Related to this, we will use $c^{\star}_i$ and $d^{\star}_i$ to denote a fixed point of a system within the presented framework, i.e.~we will use $c^{\star}_i$ and $d^{\star}_i$ to indicate a solution to the transformed stationarity conditions.

Making use of the established notation, it is straightforward to verify that the transformation specified in Eq.~\ref{eq:linTransformation}, applied to the stationarity conditions in Eqns.~\ref{eq:optim-primal2}-\ref{eq:optim-primal3} and \ref{eq:optim-dual2}-\ref{eq:optim-dual3}, can result in transformed stationarity conditions written as
\begin{align}
G_{\ell} {\underbar{c}^{\star}_{\ell}}^{(LI)} &= {\underbar{d}^{\star}_{\ell}}^{(LI)},\ \ \ell = 1,\dots, L \label{eq:transCond1}\\
m_k({\underbar{d}^{\star}_{k}}^{(CR)}) &= {\underbar{c}^{\star}_{k}}^{(CR)},\ \ k=1,\dots,K, \label{eq:transCond2}
\end{align}
where the linear map $G_{\ell}$ and the generally nonlinear map $m_k(\cdot)$ satisfy the following relationships:
\begin{align}
& \left\{ M^{(LI)}_{\ell} \left[ \begin{array}{c} \underbar{a}^{(i)}_{\ell} \\ A_{\ell} \underbar{a}^{(i)}_{\ell} \\ -A^T_{\ell}  \underbar{b}^{(o)}_{\ell} \\ \underbar{b}^{(o)}_{\ell} \end{array} \right] :  \left[ \begin{array}{c} \underbar{a}^{(i)}_{\ell} \\ \underbar{b}^{(o)}_{\ell} \end{array} \right] \in \mathbb{R}^{N^{(LI)}_{\ell}} \right\} \nonumber \\
 &\ = \left\{ \left[ \begin{array}{c} \underbar{c}^{(LI)}_{\ell}  \\ G_{\ell} \underbar{c}^{(LI)}_{\ell} \end{array} \right]: \underbar{c}^{(LI)}_{\ell} \in \mathbb{R}^{N^{(LI)}_{\ell}}\right\},\ \ell = 1, \dots, L \label{eq:Atransformation}
\end{align}
and
\begin{align}
& \left\{ M^{(CR)}_k \left[ \begin{array}{c} f_k(\underbar{y}^{(CR)}_k) \\ g_k(\underbar{y}^{(CR)}_k) \end{array} \right]: \underbar{y}^{(CR)}_k \in \mathbb{R}^{N^{(CR)}_k} \right\} \nonumber \\
& \ =  \left\{ \left[ \begin{array}{c} m_k(\underbar{d}^{(CR)}_k) \\ \underbar{d}^{(CR)}_k \end{array} \right]: \underbar{d}^{(CR)}_k \in \mathbb{R}^{N^{(CR)}_k} \right\},\ k = 1, \dots, K. \label{eq:fgTransformation}
\end{align}
Given a solution $c^{\star}_i$ and $d^{\star}_i$ to the transformed conditions written using maps in the form of Eqns.~\ref{eq:transCond1}-\ref{eq:transCond2}, the associated reduced-form primal and dual variables $a^{\star}_i$ and $b^{\star}_i$ can be obtained in a straightforward way by inverting the relationship specified by the $2\times 2$ matrices in Eq.~\ref{eq:linTransformation}. 

A significant potential obstacle in performing a change of coordinates is that for a pre-specified set of transformations $M_i$ and maps $f_k(\cdot)$, $g_k(\cdot)$ and $A_{\ell}$, there generally may not exist maps $m_k(\cdot)$ and $G_{\ell}$ that satisfy Eqns.~\ref{eq:Atransformation}-\ref{eq:fgTransformation}.  However referring to Eq.~\ref{eq:Atransformation}, there exists a class of transformations $M_i$ that will be shown in Subsection \ref{ssec:consprinc} to always result in a valid linear map $G_{\ell}$. And referring to the existence of maps $m_k(\cdot)$ satisfying Eq.~\ref{eq:fgTransformation}, a broad and useful class of generally nonlinear maps $m_k(\cdot)$ is discussed in Section \ref{sec:exarch}.

\subsection{Conservation principle}
\label{ssec:consprinc}

In designing physical systems for convex and nonconvex optimization\cite{ChuaLin}\cite{Dennis}\cite{KennedyChua}\cite{Wyatt1995} and distributed control\cite{BaranHorn}, the conservation principle resulting from Eqns.~\ref{eq:optim-primal3} and \ref{eq:optim-dual3}, specifically orthogonality between vectors of conjugate variables, is a key part of the foundation on which the systems are developed.  In electrical networks, this principle is specifically embodied by Tellegen's theorem.\cite{PenfieldSpenceDuinker}\cite{Tellegen} The conditions in Eqns.~\ref{eq:optim-primal3} and \ref{eq:optim-dual3} in particular imply
\begin{equation}
\sum_{i=1}^Na_i b_i = \sum_{\ell = 1}^L \langle \underbar{a}^{(i)}_{\ell} , -A^T_{\ell}  \underbar{b}^{(o)}_{\ell} \rangle + \langle A_{\ell} \underbar{a}^{(i)}_{\ell} , \underbar{b}^{(o)}_{\ell} \rangle \label{eq:orthoStep2} = 0.
\end{equation}

Viewing the left-hand side of Eq.~\ref{eq:orthoStep2} as a quadratic form, it can be shown to be isomorphic to the quadratic form composing the left-hand side of the following conservation principle:\cite{tbaran-phd}
\begin{equation}
\vspace{-1mm}
\label{eq:consPseudoPower}
\sum_{i=1}^Nc_i^2 - d_i^2 = 0.
\end{equation}
Eq.~\ref{eq:consPseudoPower} is similar to the statement of conservation of pseudopower in the wave-digital class of signal processing structures, and within that and other classes of systems is the foundation for analyzing stability and robustness in the presence of delay elements.\cite{Fettweis}\cite{DeprettereDewilde}\cite{RaoKailath}

Motivated by this and (R1), we specifically require that the variables $c_i$ and $d_i$ satisfy Eq.~\ref{eq:consPseudoPower}, and in particular that the $2\times 2$ matrices $M_i$ in Eq.~\ref{eq:linTransformation} be chosen so that the resulting interconnection elements $G_{\ell}$ are orthonormal matrices.  This requirement, combined with dissipation in the constitutive relations, underlies the discussion of algorithm convergence in Part II \cite{BaranLahlouPartII}.  As the stationarity conditions in Eqns.~\ref{eq:optim-primal3} and \ref{eq:optim-dual3} imply Eq.~\ref{eq:orthoStep2}, which as a quadratic form is isomorphic to Eq.~\ref{eq:consPseudoPower} using transformations of the form of Eq.~\ref{eq:linTransformation},\cite{tbaran-phd} we are ensured that such matrices $G_{\ell}$ satisfying Eq.~\ref{eq:Atransformation} will exist.

\vspace{-1mm}
\section{Example architecture elements}
\label{sec:exarch}

Figs.~\ref{fig:smallTable} and \ref{fig:bigTable} depict interconnection elements and constitutive relations that respectively satisfy Eqns.~\ref{eq:Atransformation} and \ref{eq:fgTransformation}.  A distributed, asynchronous optimization algorithm may be realized by connecting the constitutive relations in Fig.~\ref{fig:bigTable} to the interconnection elements in Fig.~\ref{fig:smallTable}  and eliminating algebraic loops as discussed previously using linear algebraic reduction and synchronous or asynchronous delays.  In Part II \cite{BaranLahlouPartII} we provide several examples of algorithms developed using this general strategy.

\FloatBarrier

\bibliographystyle{IEEEbib}
\nocite{*}

\bibliography{refs_part1}

\begin{thebibliography}{10}

\bibitem{NealBoyd}
N.~Parikh and S.~Boyd,
\newblock ``Block splitting for distributed optimization,''
\newblock {\em Mathematical Programming Computation}, 2014.

\bibitem{WeiOzdaglar1}
E.~Wei and A.~Ozdaglar,
\newblock ``Distributed alternating direction method of multipliers,''
\newblock in {\em Decision and Control (CDC), 2012 IEEE 51st Annual Conference
  on}, Dec 2012, pp. 5445--5450.

\bibitem{ForeroCanoGiannakis}
P.~A. Forero, A.~Cano, and G.~B. Giannakis,
\newblock ``Consensus-based distributed support vector machines,''
\newblock {\em J. Mach. Learn. Res.}, 2010.

\bibitem{BaranLahlouPartII}
T.~A. Baran and T.~A. Lahlou,
\newblock ``Conservative signal processing architectures for asynchronous,
  distributed optimization part {II}: Example systems,''
\newblock in {\em Proc. of IEEE Global Conference on Signal and Information
  Processing (GlobalSIP)}, 2014.

\bibitem{Willems1}
J.~C. Willems,
\newblock ``Dissipative dynamical systems part {I}: General theory,''
\newblock {\em Archive for Rational Mechanics and Analysis}, vol. 45, pp.
  321--351, jan 1972.

\bibitem{Millar}
W.~Millar,
\newblock ``Some general theorems for non-linear systems possessing
  resistance,''
\newblock {\em Philosophical Magazine Series 7}, vol. 42, no. 333, pp.
  1150--1160, 1951.

\bibitem{PenfieldSpenceDuinker}
P.~Penfield, R.~Spence, and S.~Duinker,
\newblock {\em Tellegen's Theorem and Electrical Networks},
\newblock The MIT Press, 1970.

\bibitem{ChuaLin}
L.~O. Chua and G.~N. Lin,
\newblock ``Nonlinear programming without computation,''
\newblock {\em Circuits and Systems, IEEE Transactions on}, vol. 31, no. 2, pp.
  182--188, Feb 1984.

\bibitem{Dennis}
J.~B. Dennis,
\newblock {\em Mathematical Programming and Electrical Networks},
\newblock Ph.D. thesis, Massachusetts Institute of Technology, 1958.

\bibitem{KennedyChua}
M.~P. Kennedy and L.~O. Chua,
\newblock ``Neural networks for nonlinear programming,''
\newblock {\em Circuits and Systems, IEEE Transactions on}, vol. 35, no. 5, pp.
  554--562, May 1988.

\bibitem{Wyatt1995}
J.~Wyatt,
\newblock ``Little-known properties of resistive grids that are useful in
  analog vision chip designs,''
\newblock {\em Vision Chips: Implementing Vision Algorithms with Analog VLSI
  Circuits}, pp. 72--89, 1995.

\bibitem{tbaran-phd}
T.~A. Baran,
\newblock {\em Conservation in Signal Processing Systems},
\newblock Ph.D. thesis, Massachusetts Institute of Technology, 2012.

\bibitem{Willems2}
J.~C. Willems,
\newblock ``The behavioral approach to open and interconnected systems,''
\newblock {\em Control Systems, IEEE}, vol. 27, no. 6, pp. 46--99, Dec 2007.

\bibitem{BaranHorn}
T.~A. Baran and B.~K.~P. Horn,
\newblock ``A robust signal-flow architecture for cooperative vehicle density
  control,''
\newblock in {\em Acoustics, Speech and Signal Processing (ICASSP), 2013 IEEE
  International Conference on}, May 2013, pp. 2790--2794.

\bibitem{Tellegen}
B.~D.~H. Tellegen,
\newblock ``A general network theorem, with applications,''
\newblock Tech. {R}ep., Philips Research Reports, Philips Research Reports.

\bibitem{Fettweis}
A.~Fettweis,
\newblock ``Wave digital filters: Theory and practice,''
\newblock {\em Proceedings of the IEEE}, vol. 74, no. 2, pp. 270--327, Feb
  1986.

\bibitem{DeprettereDewilde}
E.~Deprettere and P.~Dewilde,
\newblock ``Orthogonal cascade realization of real multiport digital filters,''
\newblock {\em International Journal of Circuit Theory and Applications}, vol.
  8, no. 3, pp. 245--272, 1980.

\bibitem{RaoKailath}
S.~K. Rao and T.~Kailath,
\newblock ``Orthogonal digital filters for vlsi implementation,''
\newblock {\em Circuits and Systems, IEEE Transactions on}, vol. 31, no. 11,
  pp. 933--945, Nov 1984.

\bibitem{baran2010linear}
T.~Baran, D.~Wei, and A.~V. Oppenheim,
\newblock ``Linear programming algorithms for sparse filter design,''
\newblock {\em Signal Processing, IEEE Transactions on}, vol. 58, no. 3, pp.
  1605--1617, 2010.

\bibitem{BoydNealChu}
S.~Boyd, N.~Parikh, E.~Chu, B.~Peleato, and J.~Eckstein,
\newblock ``Distributed optimization and statistical learning via the
  alternating direction method of multipliers,''
\newblock {\em Found. Trends Mach. Learn.}, 2011.

\end{thebibliography}

\end{document}